\documentclass[12pt,a4paper]{article}
\usepackage{amsmath,amssymb,theorem}

\newcommand{\frC}{\mathfrak{C}}
\newcommand{\fra}{\mathfrak{a}}
\newcommand{\frb}{\mathfrak{b}}

\newcommand{\frf}{\mathfrak{f}}
\newcommand{\Sine}{\mathcal{S}}
\newcommand{\T}{\mathcal{T}}

\newcommand{\Z}{\mathbb{Z}}
\newcommand{\Q}{\mathbb{Q}}
\newcommand{\R}{\mathbb{R}}

\newcommand{\eps}{\varepsilon}

\newcommand{\base}[1]{\langle #1\rangle}
\newcommand{\pcf}[1]{\bigl((#1)\bigr)}
\newcommand{\mcf}[1]{\bigl[[#1]\bigr]}

\DeclareMathOperator{\Cl}{\mathit{Cl}}
\DeclareMathOperator{\Li}{\mathrm{Li}}
\DeclareMathOperator{\re}{\mathrm{Re}}

\numberwithin{equation}{subsection}

\theorembodyfont{\slshape}
\newtheorem{thm}{Theorem}[subsection]
\newtheorem{prop}[thm]{Proposition}
\newtheorem{cor}[thm]{Corollary}
\newtheorem{lem}[thm]{Lemma}
\newtheorem{defn}[thm]{Definition}
\newtheorem{rem}[thm]{Remark}

\newenvironment{proof}
	{\textbf{Proof.}}
	{\hfill\rule{0.5em}{2ex}\par\bigbreak}

\begin{document}
\title
{On Kronecker limit formulas for real quadratic fields}
\author{Shuji Yamamoto}
\maketitle

\begin{abstract}
Let $\zeta(s,\frC)$  be the partial zeta function 
attached to a ray class $\frC$ of a real quadratic field. 
We study this zeta function at $s=1$ and $s=0$, 
combining some ideas and methods due to Zagier and Shintani. 
The main results are (1) a generalization of Zagier's formula 
for the constant term of the Laurent expansion at $s=1$, 
(2) some expressions for the value and the first derivative at $s=0$, 
related to the theory of continued fractions, and 
(3) a simple description of the behavior of Shintani's invariant 
$X(\frC)$, which is related to $\zeta'(0,\frC)$, 
when we change the signature of $\frC$. 
\end{abstract}

\tableofcontents

\section{Introduction}
Let $K$ be a number field and $\chi$ a ray class character. 
A Kronecker limit formula is an expression of the value 
(or the Laurent coefficient of degree $0$ when $\chi=1$) 
of the $L$-function $L(s,\chi)$ at $s=1$. 
When $K$ is the rational number field or an imaginary quadratic field, 
such formulas are classical and well-understood, with deep applications 
in number theory. The case of a real quadratic field, which we consider 
in this article, has also been studied by many authors. 
We mainly try to mix some ideas and methods of 
Zagier \cite{Zagier75,Zagier75b,Zagier77} 
and Shintani \cite{Shintani76,Shintani77,Shintani78}. 

In the following, let $K$ be real quadratic. 
For a narrow ray class $\frC$ of $K$, we denote by $\rho(\frC)$ 
the $0$-th Laurent coefficient at $s=1$ of the partial zeta function 
\[\zeta(s,\frC)=\sum_{\fra\in\frC}N(\fra)^{-s}. \]
We also call an expression of $\rho(\frC)$ a Kronecker limit formula, 
since they are essentially equivalent by the relation 
\[L(s,\chi)=\sum_{\frC\in\Cl_K(\frf)}\chi(\frC)\zeta(s,\frC),\]
where $\Cl_K(\frf)$ denotes the ray class group of modulus $\frf$ and 
$\chi$ is a character of it. 
Zagier \cite{Zagier75} proved such a formula 
when $\frC$ is a narrow ideal class, using the theory of continued fractions 
as a fundamental tool. 
Our first main result is an extension of Zagier's formula 
to narrow ray classes $\frC$ of an arbitrary modulus $\frf\subset O_K$ 
(Theorem \ref{thm:s=1}). 

On the other hand, Shintani \cite{Shintani77,Shintani77b} explained 
how to exploit the functional equation of $L$-functions 
to reduce the problem to the study of the behavior at $s=0$. 
Now let us recall it. 
Let $\chi$ be a ray class character of modulus $\frf$. 
Since $K$ is real quadratic, there are four types of signature for $\chi$: 
$\bigl(\chi(\frC_1),\chi(\frC_2)\bigr)=(\pm 1,\pm 1)$, 
where $\frC_1$ and $\frC_2$ are the ray classes defined by 
\begin{equation}\label{eq:C_i}
\begin{split}
\frC_1=[(\mu_1)],\quad \mu_1\in 1+\frf,\ \mu_1<0,\ \mu_1'>0,\\
\frC_2=[(\mu_2)],\quad \mu_2\in 1+\frf,\ \mu_2>0,\ \mu_2'<0.
\end{split}
\end{equation}
We write $b_\chi\in\{0,1,2\}$ the number of $+1$ in the signature of $\chi$, 
and put 
\[\Lambda(s,\chi)=\bigl(DN(\frf)\bigr)^{s/2}
\Gamma_\R(s)^{b_\chi}\Gamma_\R(s+1)^{2-b_\chi}L(s,\chi). \]
Here $D$ denotes the discriminant of $K$, and 
$\Gamma_\R(s)=\pi^{-s/2}\Gamma(s/2)$.
Then, if $\chi$ is a primitive character, there is the functional equation 
\begin{equation}\label{eq:FE}
\Lambda(s,\chi)=W(\chi)\Lambda(1-s,\chi^{-1}),
\qquad \bigl|W(\chi)\bigr|=1. 
\end{equation}
From this, we obtain the relation 
\begin{equation}\label{eq:s=1and0}
L(1,\chi)=C_\chi L^{(b_\chi)}(0,\chi^{-1}), \quad 
C_\chi=\frac{2^{b_\chi}\pi^{2-b_\chi}W(\chi)}{\sqrt{DN(\frf)}\ b_\chi!}.
\end{equation}
(When $\chi=1$, the left hand side must be replaced by the $0$-th Laurent 
coefficient.) Hence, if we can evaluate 
\begin{gather}
\zeta(0,\frC)-\zeta(0,\frC\frC_1)-\zeta(0,\frC\frC_2)
+\zeta(0,\frC\frC_1\frC_2),\label{eq:Intro_zeta1}\\
\zeta'(0,\frC)-\zeta'(0,\frC\frC_1)+\zeta'(0,\frC\frC_2)
-\zeta'(0,\frC\frC_1\frC_2),\label{eq:Intro_zeta2}\\
\zeta'(0,\frC)+\zeta'(0,\frC\frC_1)-\zeta'(0,\frC\frC_2)
-\zeta'(0,\frC\frC_1\frC_2),\label{eq:Intro_zeta3}\\
\zeta''(0,\frC)+\zeta''(0,\frC\frC_1)+\zeta''(0,\frC\frC_2)
+\zeta''(0,\frC\frC_1\frC_2),\label{eq:Intro_zeta4}
\end{gather}
a Kronecker limit formula is obtained for $\chi$ of signature 
$(-1,-1)$, $(-1,+1)$, $(+1,-1)$ or $(+1,+1)$, respectively. 

We deal with $\zeta(0,\frC)$ and $\zeta'(0,\frC)$ by using 
a quite general method given by Shintani \cite{Shintani76,Shintani77}. 
On the other hand, to the author's knowledge, there has been almost nothing 
known about the second or higher derivatives of $\zeta(s,\frC)$ at $s=0$, 
except for the pioneering work of Yoshida \cite[Appendix II]{Yoshida03}. 
For this reason, in this paper, we do not consider 
the case of signature $(+1,+1)$ at all. 

The method of Shintani mentioned above is based on 
a suitable choice of a cone decomposition of the first quadrant of $\R^2$. 
In the actual investigations of real quadratic fields, 
he mainly used the simplest one, which was spanned by $1$ and 
the totally positive fundamental unit (see \cite{Shintani77,Shintani78}). 
In this paper, instead, we prefer to choose one which 
is induced from the continued fraction, following Zagier, 
and obtain a generalization of his formula \cite[(3.3)]{Zagier77} 
for $\zeta(0,\frC)$. 
An advantage of this choice, aside from the interesting relation itself 
to the theory of continued fractions, 
is the possibility to compare the data associated with $\frC$ and 
$\frC\frC_2$. This was exploited by Zagier \cite[\S8]{Zagier75} 
in his proof of Meyer's theorem about $\rho(\frC)-\rho(\frC\frC_2)$ 
for narrow ideal classes $\frC$. 

The central subjects of \S\ref{sec:zeta'(0,C)} are the invariants 
\[X(\frC)=\exp\bigl(-\zeta'(0,\frC)+\zeta'(0,\frC\frC_1\frC_2)\bigr), \]
first studied by Shintani \cite{Shintani77,Shintani77b,Shintani78} 
(although our definition of $X(\frC)$ is the inverse of his). 
They have (or should have, at least) the great importance 
in the arithmetic of real quadratic fields, 
because of the Stark-Shintani conjecture which claims that 
they (or appropriate powers of them) are units of certain class fields 
over $K$ and generate them. 
Suggested by (\ref{eq:Intro_zeta2}) and (\ref{eq:Intro_zeta3}), 
we compare $X(\frC)$ and $X(\frC\frC_2)$. 
Then Zagier's cone decomposition again allows us to obtain a beautiful 
relation (Theorem \ref{thm:X_i(C^*)}), 
which leads to an expectation about contributions of 
infinite places to the value $L(1,\chi)$ (see Corollary \ref{cor:b_chi=1} 
and Remark \ref{rem:Contribution}). 

\bigbreak

The outline of the present paper is as follows. 
In \S\ref{sec:KLF}, we prove a Kronecker limit formula for ray classes 
(Theorem \ref{thm:s=1}), 
generalizing Zagier's for narrow ideal classes. 
The key point of the proof is the decomposition of the partial zeta functions 
given in Proposition \ref{prop:zetaDecomp}, 
which is also the basis of the discussions in \S\S\ref{sec:zeta(0,C)} and 
\ref{sec:zeta'(0,C)}. 

\S\ref{sec:zetas} summarizes the formulas about the values and 
the first derivatives at $s=0$ of several types of zeta functions. 
Moreover, in \ref{subsec:DoubleSine}, we give the definition and 
proofs of some elementary properties of the double sine function, 
which is fundamental in \S\ref{sec:zeta'(0,C)}. 

\S\ref{sec:zeta(0,C)} is devoted to the study of the values $\zeta(0,\frC)$, 
especially the elementary expressions of those values 
(Theorem \ref{thm:zeta(0,C)}). 
We also give the descriptions of 
the data attached to $\overline{\frC}=\frC\frC_1\frC_2$ and 
to $\frC^*=\frC\frC_2$ in \ref{subsec:CvsCbar} and \ref{subsec:CvsC^*}, 
respectively. 

In \S\ref{sec:zeta'(0,C)}, we study the invariant $X(\frC)$ 
as already mentioned. We obtain an expression of $X(\frC)$ 
by the double sine functions (Theorem \ref{thm:X(C)}), 
by using the description in \ref{subsec:CvsCbar}.  
Furthermore, we deduce a simple relation between $X(\frC)$ and $X(\frC^*)$ 
(Theorem \ref{thm:X_i(C^*)}) from Proposition \ref{prop:zz^*}. 

\bigbreak

\textbf{Acknowledgements.} 
The author would like to express his gratitude 
to his advisor T.~Tsuji for the constant encouragement and 
helpful suggestions. 
The author thanks T.~Ito and T.~Taniguchi 
for valuable comments and discussions. 

\subsection{Notation} 
Throughout the paper, $K$ denotes a real quadratic field of discriminant $D$. 
The conjugate of $x\in K$ is denoted by $x'$. 
We fix an embedding of $K$ into $\R$. 
For a subset $X$ of $K$, $X_+$ means the set of totally positive elements 
of $X$. 

For an integral ideal $\frf$ of $K$, denote the narrow ray class group 
of modulus $\frf$ by $\Cl_K(\frf)$, and let $\eps_\frf$ be the generator 
of the group $\bigl(O_K^\times\cap(1+\frf)\bigr)_+$, 
which is greater than $1$. Totally positive fundamental unit $\eps_{O_K}$ is 
simply denoted by $\eps$. 

If $x$ is a real number, we define $\langle x\rangle$ (resp. $\{x\}$) 
to be the number $t$ such that $x-t\in\Z$ and 
$0<t\leq 1$ (resp. $0\leq t<1$). This must not be confused with 
the notation $\base{a,b}$, which means the $\Z$-linear span of $a$ and $b$.

\section{A Kronecker limit formula for a ray class}\label{sec:KLF}
In this section, we prove a Kronecker limit formula similar to 
Zagier's in \cite{Zagier75}, for a ray class $\frC\in\Cl_K(\frf)$ of 
arbitrary conductor $\frf$. That is a formula for the constant term 
of the Laurent expansion at $s=1$ of the partial zeta function 
\[\zeta(s,\frC)=\sum_{\fra\in\frC,\,\fra\subset O_K}N(\fra)^{-s}.\]

\subsection{A cone decomposition}\label{subsec:ConeDecomp}
In the following, we use some results on continued fractions. 
For the proofs and further discussions on this theory, 
we refer the reader to Zagier's paper \cite{Zagier75} or his lecture note 
\cite{Zagier81}. 

We choose an integral ideal $\fra$ belonging to the class $\frC$. 
Then there exists a fractional ideal $\frb$ of the form 
\[\frb=\base{1,\omega},\quad 0<\omega'<1<\omega, \]
which is in the narrow ideal class of the ideal $\fra^{-1}\frf$, 
i.e.\ there is a totally positive number $z\in K^\times$ 
satisfying $\frb=(z)\fra^{-1}\frf$. 
Fix such $\fra$, $\frb$ and $z$. 

First, any integral ideal in $\frC$ can be written as 
$(\alpha)\fra$, where $\alpha\in K$ is totally positive and satisfies 
$\alpha-1\in\fra^{-1}\frf$. Hence 
\[\zeta(s,\frC)=\sum_{\alpha\in(1+\fra^{-1}\frf)_+/\base{\eps_\frf}}
N\bigl((\alpha)\fra\bigr)^{-s}. \]
Moreover, multiplying each $\alpha$ by $z$, we obtain 
\begin{equation}\label{eq:SumX}
\begin{split}
\zeta(s,\frC)
&=\sum_{\beta\in(z+\frb)_+/\base{\eps_\frf}}
N\bigl((\beta)\frb^{-1}\frf\bigr)^{-s}\\
&=N(\frb^{-1}\frf)^{-s}\sum_{\beta\in X}N(\beta)^{-s},
\end{split}
\end{equation}
where the last sum is taken for the set 
\[X=\{x+y\eps_\frf^{-1}\in z+\frb\mid x>0,\ y\geq 0\},\]
which is a system of representatives for $(z+\frb)_+/\base{\eps_\frf}$. 
We decompose this set by using the theory of continued fractions. 

From the condition $0<\omega'<1<\omega$, we have a purely periodic 
`minus' continued fraction expansion 
\[\omega=\mcf{b_0,\ldots,b_{m-1}}
:=b_0-\cfrac{1}{b_1-\cdots\cfrac{1}{b_{m-1}-\cfrac{1}{b_0-\cdots}}}
\qquad(b_k\geq 2). \]
We extend the sequence $\{b_k\}$ by the periodicity $b_{k+m}=b_k$ 
for all $k\in\Z$, and set $\omega_k=\mcf{b_k,\ldots,b_{k+m-1}}$. 
We also define the sequence $\{A_k\}$ by 
\[A_0=1, \quad A_{k+1}=A_k/\omega_{k+1}\quad (k\in\Z).\]
Since $\omega_k=b_k-\omega_{k+1}^{-1}$, we have 
\begin{equation}\label{eq:A_k+1}
A_{k+1}=A_k(b_k-\omega_k)=b_kA_k-A_{k-1}. 
\end{equation}
Hence 
\[\base{A_{k+1},A_k}=\base{A_k,A_{k-1}}=\cdots=\base{A_0,A_{-1}}
=\base{1,\omega_0}=\frb,\]
and there is a unique pair $(x_k,y_k)$ of rational numbers which satisfies
\[0<x_k\leq 1,\quad 0\leq y_k<1,\quad x_kA_{k-1}+y_kA_k\in z+\frb\]
for each $k\in\Z$. Since 
\[\begin{split}
x_{k+1}A_k+y_{k+1}A_{k+1}&\equiv x_kA_{k-1}+y_kA_k\\
&=(b_kx_k+y_k)A_k-x_kA_{k+1}\pmod{\frb}, 
\end{split}\]
they satisfy 
\begin{equation}\label{eq:xy_k+1}
x_{k+1}=\langle b_kx_k+y_k\rangle,\quad y_{k+1}=1-x_k. 
\end{equation}

\begin{defn}\upshape
We call the sequence $\bigl\{(\omega_k,x_k,y_k)\bigr\}$ 
the \textit{decomposition datum} associated with $\frC$. 
\end{defn}

\begin{rem}\upshape\label{rem:DDCan}
There are other candidates for the choice of $(\fra,\omega,z)$. 
For example, the all candidates of $\omega$ are $\omega_{k_0}$ 
for $k_0\in\Z$. In general, 
if we replace the choice of $(\fra,\omega,z)$ by another candidate, 
the sequence $\bigl\{(\omega_k,x_k,y_k)\bigr\}$ is replaced by 
$\bigl\{(\omega_{k+k_0},x_{k+k_0},y_{k+k_0})\bigr\}$ for some $k_0\in\Z$. 
In other words, the decomposition datum associated with $\frC$ is 
determined up to shift of the index. 
\end{rem}

\begin{lem}\label{lem:XDecomp}
We have the disjoint decomposition 
\[X=\coprod_{k=1}^{rm}
\bigl\{(x_k+p)A_{k-1}+(y_k+q)A_k\mid p,q\in\Z,\ p,q\geq 0\bigr\},\]
where $r=\log\eps_\frf/\log\eps$. 
\end{lem}
\begin{proof}
It is easy to see the decomposition 
\[(K\otimes_\Q\R)_+=\coprod_{k\in\Z}
\{pA_{k-1}+qA_k\mid p,q\in\R,\ p>0,\ q\geq 0\} \]
of the first quadrant in $K\otimes_\Q\R\cong\R^2$. Therefore 
\[(z+\frb)_+=\coprod_{k\in\Z}
\bigl\{(x_k+p)A_{k-1}+(y_k+q)A_k\mid p,q\in\Z,\ p,q\geq 0\bigr\}. \]
On the other hand, one finds the fact that $A_m=\eps^{-1}$ 
in Zagier's paper \cite[section 6]{Zagier75}. 
Hence we have $A_{rm}=\eps_\frf^{-1}$, 
and the claim follows. 
\end{proof}

\begin{prop}\label{prop:zetaDecomp}
For each $k\in\Z$, let $Q_k$ be the quadratic form defined by 
\[Q_k(x,y)=\frac{(x\omega_k+y)(x\omega_k'+y)}{\omega_k-\omega_k'}. \]
Then 
\[\zeta(s,\frC)=\bigl(D^{1/2}N(\frf)\bigr)^{-s}
\sum_{k=1}^{rm}Z_{Q_k}(s,x_k,y_k), \]
where 
\[Z_Q(s,x,y)=\sum_{p,q=0}^\infty Q(x+p,y+q)^{-s}. \]
\end{prop}
\begin{proof}
From (\ref{eq:SumX}) and Lemma \ref{lem:XDecomp}, we obtain 
\begin{align*}
\zeta(s,\frC)&=N(\frb^{-1}\frf)^{-s}\sum_{k=1}^{rm}\sum_{p,q=0}^\infty 
N\bigl((x_k+p)A_{k-1}+(y_k+q)A_k\bigr)^{-s}\\
&=N(\frb^{-1}\frf)^{-s}\sum_{k=1}^{rm}\sum_{p,q=0}^\infty 
\bigl\{N(A_k)(\omega_k-\omega_k')Q_k(x_k+p,y_k+q)\bigr\}^{-s}\\
&=N(\frb^{-1}\frf)^{-s}\sum_{k=1}^{rm}
\bigl\{N(A_k)(\omega_k-\omega_k')\bigr\}^{-s}Z_{Q_k}(s,x_k,y_k). 
\end{align*}
Hence it is sufficient to show that 
\[N(A_k)(\omega_k-\omega_k')=N(\frb)\sqrt{D}. \]
Since $\frb=\base{1,\omega_0}$, the right hand side is equal to 
$\omega_0-\omega_0'$. On the other hand, it holds that 
\begin{align*}
N(A_k)(\omega_k-\omega_k')&=N(A_{k+1}\omega_{k+1})
\bigl((b_k-\omega_{k+1}^{-1})-(b_k-\omega_{k+1}^{-1})'\bigr)\\
&=N(A_{k+1})(\omega_{k+1}-\omega_{k+1}')
\end{align*}
for any $k\in\Z$, and this common value is $\omega_0-\omega_0'$. 
\end{proof}

\subsection{The limit formula}
We prove the Kronecker limit formula:

\begin{thm}\label{thm:s=1}
The notation being the same as in \ref{subsec:ConeDecomp}, 
\[\lim_{s\to 1}\biggl(\bigl(D^{1/2}N(\frf)\bigr)^s\zeta(s,\frC)
-\frac{\log\eps_\frf}{s-1}\biggr)
=\sum_{k=1}^{rm}P(\omega_k,\omega_k',x_k,y_k),\]
where the function $P$ is defined by 
\begin{align*}
P(\omega,\omega',x,y)=
&F(\omega,x,y)-F(\omega',x,y)+\Li_2(\omega'/\omega)-\frac{\pi^2}{6}\\
&+\log(\omega/\omega')\biggl(-\psi(x)-\frac{\log(\omega-\omega')}{2}
+\frac{\log(\omega/\omega')}{4}\biggr). 
\end{align*}
Here $\psi(t)=\frac{\Gamma'(t)}{\Gamma(t)}$ is the logarithmic derivative 
of the gamma function, 
$\Li_2(t)=\sum_{n=1}^\infty\frac{t^n}{n^2}$ is 
the dilogarithm, and 
\begin{gather*}
F(\omega,x,y)=\int_0^\infty\biggl(\frac{e^{-yt}}{1-e^{-t}}-\frac{1}{t}\biggr)
f(\omega t,x)dt,\\
f(\omega,x)=-\int_\omega^\infty\frac{e^{-xu}}{1-e^{-u}}du. 
\end{gather*}
\end{thm}

By Proposition \ref{prop:zetaDecomp} and the fact that 
$\eps_\frf=\prod_{k=1}^{rm}\omega_k$, it suffices to prove:

\begin{prop}\label{prop:Z_Qs=1}
Let $\omega>\omega'$ be positive real numbers (here the prime does not mean 
the conjugate), and $Q$ the binary quadratic form defined by 
\[Q(x,y)=\frac{(x\omega+y)(x\omega'+y)}{\omega-\omega'}. \]
Then, for $x>0$ and $y\geq 0$, we have 
\[Z_Q(s,x,y)=\frac{\log(\omega/\omega')}{2}(s-1)^{-1}
+P(\omega,\omega',x,y)+O(s-1) \]
around $s=1$. 
\end{prop}
\begin{proof}
We use the method of Egami \cite{Egami86}, 
though Zagier's original method also works in this case. 

By Proposition 1 of \cite{Egami86}, we have 
\begin{align*}
&Z_Q(s,x,y)\\&=\frac{(\omega-\omega')^{1-s}}{\Gamma(s)^2}
\int_{\omega'}^\omega\bigl\{(\omega-u)(u-\omega')\bigr\}^{s-1}
\int_0^\infty t^{2s-1}\frac{e^{-yt}}{1-e^{-t}}
\frac{e^{-xut}}{1-e^{-ut}}dt\,du\\
&=\frac{(\omega-\omega')^{1-s}}{\Gamma(s)^2}\bigl(I_1(s)+I_2(s)\bigr), 
\end{align*}
where 
\begin{align*}
I_1(s)&=\int_{\omega'}^\omega\bigl\{(\omega-u)(u-\omega')\bigr\}^{s-1}
\int_0^\infty t^{2s-1}\biggl(\frac{e^{-yt}}{1-e^{-t}}-\frac{1}{t}\biggr)
\frac{e^{-xut}}{1-e^{-ut}}dt\,du, \\
I_2(s)&=\int_{\omega'}^\omega\bigl\{(\omega-u)(u-\omega')\bigr\}^{s-1}
\int_0^\infty t^{2s-2}\frac{e^{-xut}}{1-e^{-ut}}dt\,du. 
\end{align*}
The integral $I_1(s)$ is convergent when $\re(s)>\frac{1}{2}$, and 
\begin{align*}
I_1(1)&=\int_0^\infty\biggl(\frac{e^{-yt}}{1-e^{-t}}-\frac{1}{t}\biggr)
\int_{\omega'}^\omega\frac{te^{-xut}}{1-e^{-ut}}du\,dt\\
&=\int_0^\infty\biggl(\frac{e^{-yt}}{1-e^{-t}}-\frac{1}{t}\biggr)
\int_{\omega't}^{\omega t}\frac{e^{-xu}}{1-e^{-u}}du\,dt\\
&=F(\omega,x,y)-F(\omega',x,y). 
\end{align*}
On the other hand, $I_2(s)$ can be written as 
\[I_2(s)=\Gamma(2s-1)\,\zeta(2s-1,x)
\int_{\omega'}^\omega u^{1-2s}\bigl\{(\omega-u)(u-\omega')\bigr\}^{s-1}du.\]
Here $\zeta(s,x)=\sum_{n=0}^\infty(x+n)^{-s}$ is the Hurwitz zeta function. 
Hence the proposition is proved by combining the formulas 
\begin{align*}
(\omega-\omega')^{1-s}&=1-\log(\omega-\omega')\,(s-1)+O\bigl((s-1)^2\bigr),\\
\frac{\Gamma(2s-1)}{\Gamma(s)^2}&=1+O\bigl((s-1)^2\bigr),\\
\zeta(2s-1,x)&=\frac{1}{2}(s-1)^{-1}-\psi(x)+O(s-1),
\end{align*}
and 
\begin{align*}
&\int_{\omega'}^\omega u^{1-2s}\bigl\{(\omega-u)(u-\omega')\bigr\}^{s-1}du\\
&=\int_\alpha^1\bigl\{(1-u)(1-\alpha u^{-1})u^{-1}\bigr\}^{s-1}\frac{du}{u}
\qquad (\alpha:=\omega'/\omega)\\
&=-\log\alpha+\biggl(2\Li_2(\alpha)
+\frac{\log^2\alpha}{2}-\frac{\pi^2}{3}\biggr)(s-1)+O\bigl((s-1)^2\bigr). 
\end{align*}
To show the last one, use the formula 
$\int_\alpha^1\log(1-u)\frac{du}{u}=\Li_2(\alpha)-\frac{\pi^2}{6}.$
\end{proof}

\section{The formulas for certain zeta functions at $s=0$}\label{sec:zetas}
In this section, we review the formulas which describe the values 
and the first derivatives of certain types of zeta functions at $s=0$. 
Some of them include several functions 
(the double gamma functions and the double sine functions) 
introduced by Barnes \cite{Barnes99,Barnes01} 
and Shintani \cite{Shintani76,Shintani77,Shintani78}. 

\subsection{Zeta functions}
For $x>0$, let 
\[\zeta(s,x)=\sum_{p=0}^\infty(x+p)^{-s}\]
be the Hurwitz zeta function. Similarly, for $\omega>0$ and $z>0$, we define 
a function $\zeta_2(s,\omega,z)$, called Barnes' double zeta function, by 
\[\zeta_2(s,\omega,z)=\sum_{p,q=0}^\infty(z+p\omega+q)^{-s}
\quad\bigl(\re(s)>2\bigr). \]
Furthermore, for $\omega,\omega'>0$, $x>0$ and $y\geq 0$, we write 
\[\zeta\bigl(s,(\omega,\omega'),(x,y)\bigr)
=\sum_{p,q=0}^\infty\bigl((z+p\omega+q)(z'+p\omega'+q)\bigr)^{-s}
\ \bigl(\re(s)>1\bigr), \]
where $z=x\omega+y$ and $z'=x\omega'+y$. 
(In this section, the prime does not mean the conjugate.) 

These zeta functions are known to be meromorphically continued to 
the whole $s$-plane, and holomorphic at $s=0$. 
This fact and the following proposition are the special cases 
of Corollary to Proposition 1 of \cite{Shintani76}. 
(For (3), see also \cite{Egami86}.) 

\begin{prop}\label{prop:zeta(0)}
For $\omega,\omega'>0$, $x>0$ and $y\geq 0$, 
set $z=x\omega+y$, $z'=x\omega'+y$. Then: 
\begin{enumerate}
\item $\zeta(0,x)=-B_1(x)$. 
\item $\displaystyle \zeta_2(0,\omega,x\omega+y)
=\frac{\omega}{2}B_2(x)+B_1(x)B_1(y)+\frac{1}{2\omega}B_2(y)$. 
\item $\displaystyle \zeta\bigl(0,(\omega,\omega'),(x,y)\bigr)
=\frac{1}{2}\bigl\{\zeta_2(0,\omega,z)+\zeta_2(0,\omega',z')\bigr\}$. 
\end{enumerate}
Here $B_1(x)=x-\frac{1}{2}$ and $B_2(x)=x^2-x+\frac{1}{6}$ 
denote the first and second Bernoulli polynomials. 
\end{prop}

\subsection{The derivatives at $s=0$}
The first derivative of $\zeta(s,x)$ at $s=0$ is 
expressed by Lerch's formula 
\[\zeta'(0,x)=\log\frac{\Gamma(x)}{\sqrt{2\pi}}. \]
We define a function $G(\omega,z)$ to be the similar derivative 
for $\zeta_2$: 
\[G(\omega,z):=\zeta_2'(0,\omega,z)\quad (\omega,z>0). \]
A suitable normalization of $\exp G(\omega,z)$ is called Barnes' 
double gamma function. Basic properties of this function, 
including the fact that $\exp\bigl(-G(\omega,z)\bigr)$ can be 
continued to an entire function of $z$, were investigated by 
Barnes \cite{Barnes99,Barnes01}. 

The analogous derivative of 
$\zeta\bigl(s,(\omega,\omega'),(x,y)\bigr)$ can be 
expressed as follows: 

\begin{prop}\label{prop:zeta'(0)}
For $\omega,\omega'>0$, $x>0$ and $y\geq 0$, 
\[\zeta'\bigl(0,(\omega,\omega'),(x,y)\bigr)
=G(\omega,z)+G(\omega',z')+\frac{\omega-\omega'}{4\omega\omega'}
\log\biggl(\frac{\omega'}{\omega}\biggr)B_2(y), \]
where $z=x\omega+y$ and $z'=x\omega'+y$. 
\end{prop}

For the proof, see Proposition 3 of \cite{Shintani77}, or \cite{Egami86}. 

\subsection{The double sine function}\label{subsec:DoubleSine}
The double sine function $\Sine(\omega,z)$ is defined by 
\[\Sine(\omega,z)=\exp\bigl(G(\omega,1+\omega-z)-G(\omega,z)\bigr). \]
This function was originally introduced by Shintani 
\cite{Shintani77,Shintani78}, and recently studied by 
Kurokawa-Koyama \cite{Kurokawa-Koyama03}. 

In the following proposition, we collect several properties of 
$\Sine(\omega,z)$ which are needed later. 

\begin{prop}\label{prop:Sine}
\begin{enumerate}
\item $\Sine(\omega,1+\omega-z)=\Sine(\omega,z)^{-1}$. 
\item $\Sine(\omega,1)=\omega^{1/2}$, $\Sine(\omega,\omega)=\omega^{-1/2}$. 
\item $\Sine(\omega,z)=\Sine(1/\omega,z/\omega)$. 
\item $\Sine(\omega,z)=2\sin(\pi z)\Sine(\omega,z+\omega)
=2\sin(\pi z/\omega)\Sine(\omega,z+1)$. 
\item If $\omega>1$, 
$\displaystyle \Sine(\omega,z)=
2\sin(\pi z/\omega)\frac{\Sine(\omega-1,z)}{\Sine(1-1/\omega,z/\omega)}$.
\item If $\omega<1$, 
$\displaystyle \Sine(\omega,z)=
2\sin(\pi z)\frac{\Sine(1/\omega-1,z/\omega)}{\Sine(1-\omega,z)}$. 
\end{enumerate}
\end{prop}
\begin{proof}
(1) is clear from the definition. 

For (2), we compute as follows: 
\begin{align*}
&G(\omega,\omega)-G(\omega,1)\\
&=\frac{\partial}{\partial s}
\Biggl(\sum_{p=1}^\infty\sum_{q=0}^\infty(p\omega+q)^{-s}
-\sum_{p=0}^\infty\sum_{q=1}^\infty(p\omega+q)^{-s}\Biggr)\Biggm|_{s=0}\\
&=\frac{\partial}{\partial s}(\omega^{-s}-1)\zeta(s)\bigm|_{s=0}\\
&=-\zeta(0)\log\omega. 
\end{align*}
This leads to (2), since $\zeta(0)=-\frac{1}{2}$. 

Next, we differentiate the evident identity 
\[\zeta_2(s,1/\omega,z/\omega)=\omega^s\zeta_2(s,\omega,z)\]
to obtain 
\[G(1/\omega,z/\omega)=G(\omega,z)+\zeta_2(0,\omega,z)\log\omega. \]
Changing $z$ to $1+\omega-z$, and subtracting, we get an equality 
\[\log\frac{\Sine(1/\omega,z/\omega)}{\Sine(\omega,z)}
=\bigl\{\zeta_2(0,\omega,1+\omega-z)-\zeta_2(0,\omega,z)\bigr\}\log\omega, \]
whose right hand side vanishes by Proposition \ref{prop:zeta(0)} (2). 
This proves (3). 

To show (4), we start with another identity 
\[\zeta_2(s,\omega,z)=\zeta_2(s,\omega,z+\omega)+\zeta(s,z), \] 
which is again immediate from the definition. Then, in a similar manner 
to the proof of (3) above, we obtain the first equality of (4). 
The second one can be proved in the same way, or by combining 
the first one and (3). 

Finally, the proofs of (5) and (6) are given by beginning with 
\begin{align*}
&\zeta_2(s,\omega-1,z)\\
&=\sum_{p=0}^\infty\sum_{q=-\infty}^\infty(z+p\omega+q)^{-s}\\
&=\sum_{p=0}^\infty\sum_{q=1}^\infty(z+p\omega+q)^{-s}
+\sum_{p,q=0}^\infty\bigl(z+(p+q)\omega-q\bigr)^{-s}\\
&=\zeta_2(s,\omega,z)-\omega^{-s}\zeta(s,z/\omega)
+\omega^{-s}\zeta_2(s,1-1/\omega,z/\omega), 
\end{align*}
and repeating the method above. 
\end{proof}

\section{Formulas for $\zeta(0,\frC)$}\label{sec:zeta(0,C)}
In this section, we compute the values $\zeta(0,\frC)$ by combining 
Proposition \ref{prop:zetaDecomp} and Proposition \ref{prop:zeta(0)}, 
following the general method of Shintani \cite{Shintani76}. 
Then our special choice of the cone decomposition 
based on the continued fractions leads to a particularly simple expression, 
and allows us to analyze the multiplication by $\frC_1$ and $\frC_2$. 

We use the notation introduced in \ref{subsec:ConeDecomp}. 

\subsection{An elementary expression of $\zeta(0,\frC)$}
Here we prove the following formula, which is a generalization 
of Zagier's \cite[(3.3)]{Zagier77}. 

\begin{thm}\label{thm:zeta(0,C)}
\[\zeta(0,\frC)=\sum_{k=1}^{rm}
\biggl\{B_1(x_k)B_1(y_k)+\frac{b_k}{2}B_2(x_k)\biggr\}. \]
\end{thm}
\begin{proof}
By (2) and (3) of Proposition \ref{prop:zeta(0)}, we have 
\begin{equation}\label{eq:Z(0)}
\begin{split}
&Z_{Q_k}(0,x_k,y_k)\\
&=\frac{\omega_k+\omega_k'}{4}B_2(x_k)+B_1(x_k)B_1(y_k)+\frac{1}{4}
\biggl(\frac{1}{\omega_k}+\frac{1}{\omega_k'}\biggr)B_2(y_k). 
\end{split}
\end{equation}
Therefore, we complete the proof when we substitute this into 
Proposition \ref{prop:zetaDecomp} and compute as 
\begin{align*}
&\sum_{k=1}^{rm}\Biggl\{(\omega_k+\omega_k')B_2(x_k)+\biggl(
\frac{1}{\omega_k}+\frac{1}{\omega_k'}\biggr)B_2(y_k)\Biggr\}\\
&=\sum_{k=1}^{rm}\Biggl\{(\omega_k+\omega_k')B_2(x_k)+\biggl(
\frac{1}{\omega_{k+1}}+\frac{1}{\omega_{k+1}'}\biggr)B_2(1-y_{k+1})\Biggr\}\\
&=\sum_{k=1}^{rm}\Biggl\{\omega_k+\omega_k'+
\frac{1}{\omega_{k+1}}+\frac{1}{\omega_{k+1}'}\Biggr\}B_2(x_k)\\
&=\sum_{k=1}^{rm}2b_kB_2(x_k), 
\end{align*}
using the periodicities and recurrence relations of 
$\omega_k$ and $(x_k,y_k)$, and the identity $B_2(x)=B_2(1-x)$. 
\end{proof}

\begin{rem}\upshape
Meyer \cite{Meyer57} and Siegel \cite{Siegel80} 
proved a formula similar to the above theorem 
(see the theorems 12 and 13 in \cite[\S2]{Siegel80}). 
Their proof were based on an integral formula due to Hecke, 
and completely different from ours. 
\end{rem}

\subsection{$\frC$ versus $\overline{\frC}$}\label{subsec:CvsCbar}
Let $\frC_1$ and $\frC_2$ be the ray classes defined in (\ref{eq:C_i}). 
In the following, we write $\overline{\frC}=\frC\frC_1\frC_2$ for brevity. 
We want to compare $\zeta(0,\frC)$ and $\zeta(0,\overline{\frC})$. 
For this purpose, we need the decomposition datum 
associated with $\overline{\frC}$. 

Since $\overline{\frC}$ and $\frC$ are in the common narrow ideal class, 
we may use the same $\frb$, $\omega_k$, and $A_k$. 
On the other hand, if we choose an element $\nu\in 1+\frf$ which is 
totally negative, the ideal $(\nu)\fra$ is a representative of 
$\overline{\frC}$. (Here $\fra$ is a representative of $\frC$ fixed in 
\ref{subsec:ConeDecomp}.) 
Then, since $\frb=(z)\fra^{-1}\frf=(-z\nu)\bigl((\nu)\fra\bigr)^{-1}\frf$, 
we can take the set 
\[-z\nu+\frb=-z+(1-\nu)z+\frb=-z+\frb\]
as the counterpart of $z+\frb$, 
and hence the counterparts of $x_k$ and $y_k$ becomes 
$\langle -x_k\rangle$ and $\{-y_k\}$, respectively. 
Substituting these data to Proposition \ref{prop:zetaDecomp} and 
Theorem \ref{thm:zeta(0,C)}, we obtain 
\begin{gather}
\zeta(s,\overline{\frC})=\bigl(D^{1/2}N(\frf)\bigr)^{-s}
\sum_{k=1}^{rm}Z_{Q_k}\bigl(s,\langle -x_k\rangle,\{-y_k\}\bigr),
\label{eq:zeta(Cbar)Decomp}\\
\zeta(0,\overline{\frC})=\sum_{k=1}^{rm}
\biggl\{B_1\bigl(\langle -x_k\rangle\bigr)B_1\bigl(\{-y_k\}\bigr)
+\frac{b_k}{2}B_2\bigl(\langle -x_k\rangle\bigr)\biggr\}. 
\label{eq:zeta(0,Cbar)}
\end{gather}

In fact, we have the following: 

\begin{prop}\label{prop:zeta(0,Cbar)}
Writing $\overline{\frC}=\frC\frC_1\frC_2$, we have 
\[\zeta(0,\frC)=\zeta(0,\overline{\frC}). \]
\end{prop}
\begin{proof}
We compare expressions in Theorem \ref{thm:zeta(0,C)} 
and (\ref{eq:zeta(0,Cbar)}). Since $\langle -x_k\rangle$ is 
equal to $1-x_k$ (when $x_k\in(0,1)$) or to $x_k$ (when $x_k=1$), 
$B_2\bigl(\langle -x_k\rangle\bigr)$ is equal to 
$B_2(x_k)$ for each $k$. To deal with the terms $B_1(x_k)B_1(y_k)$, 
it is necessary to discuss some cases separately. 

First we treat the case in which $(x_k,y_k)=(1,0)$ holds for some $k$. 
This means $z\in\frb$, which happens if and only if $\frf=O_K$. 
The theorem itself is trivial in this case since $\frC=\overline{\frC}$. 

Next, we assume that $x_k$ and $y_k$ are both in the open interval 
$(0,1)$, for an index $k$. 
Then $\langle -x_k\rangle=1-x_k$ and $\{-y_k\}=1-y_k$, and hence 
\[B_1(x_k)B_1(y_k)
=B_1\bigl(\langle -x_k\rangle\bigr)B_1\bigl(\{-y_k\}\bigr). \]

Finally, from the recurrence relation (\ref{eq:xy_k+1}), 
we see that $x_k=1$ if and only if $y_{k+1}=0$, and then $x_{k+1}=y_k$ 
(note that we exclude the case $(x_k,y_k)=(1,0)$). 
In this case, we have 
\begin{align*}
&B_1(x_k)B_1(y_k)+B_1(x_{k+1})B_1(y_{k+1})\\
&=B_1\bigl(\langle -x_k\rangle\bigr)B_1\bigl(\{-y_k\}\bigr)
+B_1\bigl(\langle -x_{k+1}\rangle\bigr)B_1\bigl(\{-y_{k+1}\}\bigr),
\end{align*}
as desired. 
\end{proof}

\begin{rem}\label{rem:zeta(0,Cbar)}\upshape
Proposition \ref{prop:zeta(0,Cbar)} is also deduced from the fact that 
$L(0,\chi)=0$ for any ray class character $\chi$ 
of signature $(+1,-1)$ or $(-1,+1)$. 
\end{rem}

\subsection{$\frC$ versus $\frC^*$}\label{subsec:CvsC^*}
Let us write $\frC^*=\frC\frC_2$, and consider the relation 
between $\zeta(0,\frC)$ and $\zeta(0,\frC^*)$. 
Then we can use the fact that $L(0,\chi)=0$ 
whenever $\chi$ has the signature $(+1,+1)$ or $(-1,+1)$, to deduce that 
\[\zeta(0,\frC)+\zeta(0,\frC^*)=0\]
(compare with Remark \ref{rem:zeta(0,Cbar)}). 
Unfortunately, it seems difficult to obtain this relation 
by a direct computation as in the proof of 
Proposition \ref{prop:zeta(0,Cbar)}, 
except for the case of modulus $O_K$ which was treated by 
Zagier \cite{Zagier75b,Zagier77} (see Remark \ref{rem:zeta(0,C^*)} below). 
Here we only describe the relation between the decomposition data 
associated with $\frC$ and $\frC^*$. 
We refer the reader again to \cite{Zagier75,Zagier81} 
for the theory of continued fractions used below. 

We can assume that $\omega_0>2$, by shifting the index if necessary 
(see Remark \ref{rem:DDCan}). 
Then, putting $\xi=\omega_0-1$, we have a `plus' continued fraction 
\[\xi=\pcf{a_0,\ldots,a_{2l-1}}
:=a_0+\cfrac{1}{a_1+\cdots\cfrac{1}{a_{2l-1}+\cfrac{1}{a_0+\cdots}}}
\quad(a_j\geq 1). \]
Here $2l$ denotes the smallest even period ($l$ may be the smallest period 
if $l$ is odd). We define $a_j$ for all $j\in\Z$ by $a_{j+2l}=a_j$. 
Then the sequence $\{b_k\}$ is determined by $\{a_j\}$ as 
\begin{equation}\label{eq:a_jb_k}
b_{S_j}=a_{2j}+2,\qquad b_k=2\quad (S_j<k<S_{j+1}), 
\end{equation}
using the sequence $\{S_j\}$ defined by 
\[S_0=0,\qquad S_j=S_{j-1}+a_{2j-1}. \]
Moreover, if we set $\xi_j=\pcf{a_j,\ldots,a_{j+2l-1}}$, we have 
\begin{gather}
a_{2j}+2-\frac{1}{\omega_{S_j+1}}=\omega_{S_j}
=\xi_{2j}+1=a_{2j}+1+\frac{1}{\xi_{2j+1}},\label{eq:omegaxi}\\
\omega_k=2-\frac{1}{\omega_{k+1}}\quad(S_j<k<S_{j+1}). \label{eq:omegaRec}
\end{gather}

On the other hand, setting $\omega^*=\xi_1+1$, 
we can take $\frb^*=\base{1,\omega^*}$ as the counterpart of $\frb$ 
for $\frC^*$, since $\xi_1>1$, $-1<\xi_1'<0$, and 
\[\frb^*=\xi_1\base{1,1/\xi_1}=\xi_1\base{1,\xi_0}=\xi_1\frb. \]
For $\omega^*$, we have a continued fraction expansion 
\[\omega^*=\mcf{c_0,\ldots,c_{n-1}}\quad (c_k\geq 2). \]
The sequence $\{c_k\}$ is determined by 
\begin{gather*}
T_0=0,\qquad T_j=T_{j-1}+a_{2j}, \\
c_{T_j}=a_{2j+1}+2,\qquad c_k=2\quad (T_j<k<T_{j+1}). 
\end{gather*}
We set $\omega^*_k=\mcf{c_k,\ldots,c_{k+n-1}}$. 
The identities similar to (\ref{eq:omegaxi}) and (\ref{eq:omegaRec}) are 
\begin{gather}
a_{2j+1}+2-\frac{1}{\omega^*_{T_j+1}}=\omega^*_{T_j}
=\xi_{2j+1}+1=a_{2j+1}+1+\frac{1}{\xi_{2j+2}},\label{eq:omega^*xi}\\
\omega^*_k=2-\frac{1}{\omega^*_{k+1}}\quad(T_j<k<T_{j+1}).\label{eq:omega^*Rec}
\end{gather}

Next, let us look at the counterparts of $x_k$ and $y_k$. 
We take $(\mu_2)\fra$ as a representative of $\frC^*$, 
where $\mu_2$ is a number as in (\ref{eq:C_i}). Then, since 
\[\frb^*=(\xi_1)\frb=(z^*)\bigl((\mu_2)\fra\bigr)^{-1}\frf,
\quad z^*:=z\xi_1\mu_2\gg 0, \]
we can determine rational numbers $x^*_k$ and $y^*_k$ by 
\[x^*_kA^*_{k-1}+y^*_kA^*_k\in z^*+\frb^*=\xi_1(z+\frb). \]
Here $A^*_k$ is defined from $\omega^*_k$ in the same way as $A_k$. 

Put $z_k=x_k\omega_k+y_k$ and $z^*_k=x^*_k\omega^*_k+y^*_k$. 
We want some relation between the sequences $\{z_k\}$ and 
$\{z^*_k\}$, but any one-to-one correspondence is impossible, 
since the periods $rm$ and $rn$ of them are different in general. 
There is, however, such a relation between $\{z_{S_j}\}$ and 
$\{z^*_{T_j}\}$. 

\begin{prop}\label{prop:zz^*}
For each $j$, there are congruences 
\begin{align*}
z_{S_j}\equiv\xi_{2j}z^*_{T_{j-1}}&\mod \base{1,\xi_{2j}}, \\
\xi_{2j+1}z_{S_j}\equiv z^*_{T_j}&\mod \base{1,\xi_{2j+1}}, 
\end{align*}
\end{prop}
\begin{proof}
First, we note that 
\[A_{S_{j-1}}/A_{S_{j-1}+a}=a(\omega_{S_{j-1}+a}-1)+1\]
if $1\leq a\leq a_{2j-1}=S_j-S_{j-1}$. 
This can be verified by induction on $a$, using (\ref{eq:omegaRec}). 
When $a=a_{2j-1}$, this becomes 
\begin{align*}
A_{S_{j-1}}/A_{S_j}&=a_{2j-1}(\omega_{S_j}-1)+1=a_{2j-1}\xi_{2j}+1\\
&=\xi_{2j}\xi_{2j-1}. 
\end{align*}
In a similar way, we also obtain 
\[A^*_{T_{j-1}}/A^*_{T_j}=\xi_{2j+1}\xi_{2j}. \]
From these, we see 
\[A_{S_j}^{-1}=\xi_{2j}\xi_{2j-1}\cdots\xi_1,\quad 
A^{*\,-1}_{T_j}=\xi_{2j+1}\xi_{2j}\cdots\xi_2, \]
and hence 
\[A_{S_j}^{-1}(z+\frb)=\xi_{2j}A^{*\,-1}_{T_{j-1}}(z^*+\frb^*)
=\xi_{2j+1}^{-1}A^{*\,-1}_{T_j}(z^*+\frb^*). \]
Since the three numbers $z_{S_j}$, $\xi_{2j}z^*_{T_{j-1}}$ and 
$\xi_{2j+1}^{-1}z^*_{T_j}$ belong to this common set, 
they are congruent modulo the ideal 
\[A_{S_j}^{-1}\frb=\base{1,\omega_{S_j}}=\base{1,\xi_{2j}}
=\xi_{2j+1}^{-1}\base{1,\xi_{2j+1}}. \]
This leads to the desired congruences. 
\end{proof}

\begin{rem}\label{rem:zeta(0,C^*)}\upshape
If $\frf=O_K$, Theorem \ref{thm:zeta(0,C)} can be written as 
\[\zeta(0,\frC)=\frac{1}{12}\sum_{k=1}^m(b_k-3)
=-\frac{1}{12}m+\frac{1}{12}\sum_{j=1}^la_{2j}=\frac{1}{12}(n-m). \]
This leads, in particular, to the identity $\zeta(0,\frC)=-\zeta(0,\frC^*)$. 
\end{rem}

\section{Formulas for $\zeta'(0,\frC)$}\label{sec:zeta'(0,C)}
We keep the notation in the previous section. 

Here we consider the derivative $\zeta'(0,\frC)$, 
or rather $\zeta'(0,\frC)-\zeta'(0,\overline{\frC})$. 
We follow the method of Shintani \cite{Shintani77} in principle, 
but we can compare those values for $\frC$ and $\frC^*$ 
by virtue of the continued fraction theory. 

\subsection{The invariant $X(\frC)$}
As mentioned in the introduction, we define an invariant $X(\frC)$ 
of a ray class $\frC\in\Cl_K(\frf)$ by 
\[X(\frC):=\exp\bigl(-\zeta'(0,\frC)+\zeta'(0,\overline{\frC})\bigr). \]

\begin{thm}\label{thm:X(C)}
\[X(\frC)=\prod_{k=1}^{rm}\Sine(\omega_k,z_k)\Sine(\omega_k',z_k'). \]
\end{thm}
\begin{proof}
We combine Proposition \ref{prop:zetaDecomp} and 
Proposition \ref{prop:zeta'(0)} to find 
\begin{align*}
\zeta'(0,\frC)=\sum_{k=1}^{rm}\Biggl\{&G(\omega_k,z_k)+G(\omega_k',z_k')
+\frac{\omega_k-\omega_k'}{4\omega_k\omega_k'}
\log\biggl(\frac{\omega_k'}{\omega_k}\biggr)B_2(y_k)\\
&-\log\biggl(\frac{D^{1/2}N(\frf)}{\omega_k-\omega_k'}\biggr)
Z_{Q_k}(0,x_k,y_k)\Biggr\}. 
\end{align*}
Note that $Z_{Q_k}(0,x_k,y_k)$ was already given in (\ref{eq:Z(0)}). 
Beginning with the expression (\ref{eq:zeta(Cbar)Decomp}), 
we obtain a similar formula for $\zeta'(0,\overline{\frC})$. 
Then we can compute the difference of them by case-by-case discussion, 
exactly the same as in the proof of Proposition \ref{prop:zeta(0,Cbar)}. 

When $z\in\frb$, the theorem itself becomes almost trivial: 
since $z_k=\omega_k$ for all $k$, the right hand side can be computed as 
\[\prod_{k=1}^{m}\Sine(\omega_k,\omega_k)\Sine(\omega_k',\omega_k')
\prod_{k=1}^{m}(\omega_k\omega_k')^{-1/2}=(\eps\eps')^{-1/2}=1,\]
by Proposition \ref{prop:Sine} (2). 

Next we assume $z\notin\frb$, and consider the difference for each $k$. 
If $x_k$ and $y_k$ are both in $(0,1)$, then $\langle -x_k\rangle=1-x_k$ and 
$\{-y_k\}=1-y_k$. Hence all terms of Bernoulli polynomials are cancelled out, 
and the difference becomes 
\[G(\omega_k,1+\omega_k-z_k)-G(\omega_k,z_k)
+G(\omega_k',1+\omega_k'-z_k')-G(\omega_k',z_k') \]
as desired. 

The most subtle is the remaining case, in which $x_k=1$, $y_{k+1}=0$ and 
$x_{k+1}=y_k\in(0,1)$. We have to consider $k$ and $k+1$ simultaneously, 
and use the formula 
\begin{align*}
&G(\omega_k,\omega_k+1-y_k)+G\bigl(\omega_{k+1},(1-x_{k+1})\omega_{k+1}\bigr)\\
&=G(\omega_k,1-y_k)+G\bigl(\omega_{k+1},(1-x_{k+1})\omega_{k+1}+1\bigr)
+B_1(y_k)\log\omega_{k+1}
\end{align*}
with its conjugate. The term 
$B_1(y_k)\log(\omega_{k+1}\omega_{k+1}')$ will be cancelled with 
the $B_1$ terms from $\Z_{Q_k}(0,x_k,y_k)$ etc., since 
\[\omega_k-\omega_k'=\biggl(b_k-\frac{1}{\omega_{k+1}}\biggr)-
\biggl(b_k-\frac{1}{\omega_{k+1}'}\biggr)
=\frac{1}{\omega_{k+1}\omega_{k+1}'}(\omega_{k+1}-\omega_{k+1}'). \]
We omit the detailed computation. 
\end{proof}

\subsection{$X(\frC)$ and $X(\frC^*)$}
By Theorem \ref{thm:X(C)}, we may split $X(\frC)$ as 
\begin{equation}\label{eq:X=X_1X_2}
X(\frC)=X_1(\frC)X_2(\frC),
\end{equation}
where 
\[X_1(\frC)=\prod_{k=1}^{rm}\Sine(\omega_k,z_k),\quad 
X_2(\frC)=\prod_{k=1}^{rm}\Sine(\omega_k',z_k'). \]
These are invariants of $\frC$ and independent of the choices 
made in \ref{subsec:ConeDecomp} (see Remark \ref{rem:DDCan}). 

We want to express $X_1(\frC)$ and $X_2(\frC)$ 
by $\{\xi_j\}$ instead of $\{\omega_k\}$, 
to compare $X_i(\frC)$ and $X_i(\frC^*)$ ($i=1,2$) 
using the relations explained in \ref{subsec:CvsC^*}. 
For this purpose, it is convenient to introduce some auxiliary functions. 

For a positive irrational number $\omega$, 
and and a number of the form $z=x\omega+y$ with $x,\,y\in\Q$, we define 
\[\T_1(\omega,z)
=\Sine\bigl(\omega,\langle x\rangle\omega+\langle y\rangle\bigr),\quad
\T_2(\omega,z)=\Sine\bigl(\omega,\{x\}\omega+\langle y\rangle\bigr). \]

\begin{lem}\label{lem:T_i}
Let $\omega>0$ be an irrational number, and $x\in(0,1]$ and $y\in[0,1)$ 
rational numbers. Put $z=x\omega+y$. 
\begin{enumerate}
\item If $\omega>1$, 
\[\Sine(\omega,z)=\frac{\T_1(\omega-1,z)}{\T_1(1-1/\omega,z/\omega)}.\]
\item If $\omega<1$, 
\[\Sine(\omega,z)=\frac{\T_2(1/\omega-1,z/\omega)}{\T_2(1-\omega,z)}.\]
\end{enumerate}
\end{lem}
\begin{proof}
If $\omega>1$, we combine (4) and (5) of Proposition \ref{prop:Sine} 
to obtain 
\[\Sine(\omega,z)=\frac{\Sine(\omega-1,z)}
{\Sine\Bigl(1-\frac{1}{\omega},\frac{z}{\omega}+1-\frac{1}{\omega}\Bigr)}
=\frac{\Sine(\omega-1,z-1)}
{\Sine\Bigl(1-\frac{1}{\omega},\frac{z}{\omega}-\frac{1}{\omega}\Bigr)}. \]
This leads to (1), since 
\[z=x(\omega-1)+x+y,\quad \frac{z}{\omega}+1-\frac{1}{\omega}
=(1-y)\biggl(1-\frac{1}{\omega}\biggr)+x+y. \]
(2) can be proved in a similar way. 
\end{proof}

\begin{prop}\label{prop:X_i(C)}
$X_1(\frC)$ and $X_2(\frC)$ can be written as 
\[X_1(\frC)=\prod_{j=1}^{rl}\frac{\T_1\bigl(\xi_{2j},z_{S_j}\bigr)}
{\T_1\bigl(1/\xi_{2j+1},z_{S_j}\bigr)}, \quad
X_2(\frC)=\prod_{j=1}^{rl}
\frac{\T_2\bigl(-1/\xi_{2j+1}',z_{S_j}'\bigr)}
{\T_2\bigl(-\xi_{2j}',z_{S_j}'\bigr)}. \]
\end{prop}
\begin{proof}
Since 
\begin{gather*}
z_k=x_k\omega_k+y_k=A_k^{-1}(x_kA_{k-1}+y_kA_k)\in A_k^{-1}(z+\frb), \\
z_{k+1}/\omega_{k+1}\in\omega_{k+1}^{-1}A_{k+1}^{-1}(z+\frb)
=A_k^{-1}(z+\frb), 
\end{gather*}
we have 
\[z_k\equiv z_{k+1}/\omega_{k+1}\bmod A_k^{-1}\frb=A_k^{-1}\base{A_k,A_{k-1}}
=\base{1,\omega_k}=\base{1,1-1/\omega_{k+1}}. \]
Hence, by using \ref{lem:T_i} (1), we obtain 
\begin{align*}
X_1(\frC)&=\prod_{k=1}^{rm}
\frac{\T_1(\omega_k-1,z_k)}{\T_1(1-1/\omega_k,z_k/\omega_k)}\\
&=\prod_{k=1}^{rm}
\frac{\T_1(\omega_k-1,z_k)}{\T_1(1-1/\omega_{k+1},z_{k+1}/\omega_{k+1})}\\
&=\prod_{k=1}^{rm}
\frac{\T_1(\omega_k-1,z_k)}{\T_1(1-1/\omega_{k+1},z_k)}. 
\end{align*}
Now we can prove the first formula of the theorem by substituting 
\begin{gather*}
\omega_{S_j}-1=\xi_{2j},\qquad 1-1/\omega_{S_j+1}=1/\xi_{2j+1},\\
\omega_k-1=1-1/\omega_{k+1}\qquad(S_j<k<S_{j+1}),
\end{gather*}
which are deduced from (\ref{eq:omegaxi}) and (\ref{eq:omegaRec}). 
The second one can be proved in a similar manner. 
\end{proof}

\begin{thm}\label{thm:X_i(C^*)}
Writing $\frC^*=\frC\frC_2$, we have 
\[X_1(\frC)=X_1(\frC^*)^{-1},\qquad X_2(\frC)=X_2(\frC^*). \]
\end{thm}
\begin{proof}
We may apply Proposition \ref{prop:X_i(C)} to $\frC^*$ to obtain 
\begin{equation*}
X_1(\frC^*)
=\prod_{j=1}^{rl}\frac{\T_1\bigl(\xi_{2j+1},z^*_{T_j}\bigr)}
{\T_1\bigl(1/\xi_{2j},z^*_{T_{j-1}}\bigr)}, \quad
X_2(\frC^*)
=\prod_{j=1}^{rl}
\frac{\T_2\bigl(-1/\xi_{2j}',(z^*_{T_{j-1}})'\bigr)}
{\T_2\bigl(-\xi_{2j+1}',(z^*_{T_j})'\bigr)}. 
\end{equation*}
Here we also use the periodicity 
\[\xi_{2(j+l)}=\xi_{2j},\qquad \omega^*_{T_{j+rl}}=\omega^*_{T_j+rn}
=\omega^*_{T_j}. \]
Hence it is sufficient to prove 
\begin{align*}
\T_1\bigl(\xi_{2j},z_{S_j}\bigr)&=\T_1\bigl(1/\xi_{2j},z^*_{T_{j-1}}\bigr),\\
\T_1\bigl(1/\xi_{2j+1},z_{S_j}\bigr)&=\T_1\bigl(\xi_{2j+1},z^*_{T_j}\bigr),\\
\T_2\bigl(-\xi_{2j}',z_{S_j}'\bigr)
&=\T_2\bigl(-1/\xi_{2j}',(z^*_{T_{j-1}})'\bigr)^{-1},\\
\T_2\bigl(-1/\xi_{2j+1}',z_{S_j}'\bigr)
&=\T_2\bigl(-\xi_{2j+1}',(z^*_{T_j})'\bigr)^{-1}. 
\end{align*}
In view of (1) and (3) of Proposition \ref{prop:Sine}, 
these follow from the congruences in Proposition \ref{prop:zz^*}. 
\end{proof}

\begin{rem}\upshape
\begin{enumerate}
\item By interchanging the role of two infinite places, we have 
\[X_1(\frC)=X_1(\frC\frC_1),\qquad X_2(\frC)=X_2(\frC\frC_1)^{-1}. \]
\item Yoshida showed the above formula (or, at least, a similar result) 
under some special assumptions 
(see the proposition 6.2 in \cite[Chap.~III]{Yoshida03}). 
\end{enumerate}
\end{rem}

By using Theorem \ref{thm:X_i(C^*)}, we immediately obtain 
an expression of $L(1,\chi)$ when $b_\chi=1$, 
as indicated in the introduction. 

\begin{cor}\label{cor:b_chi=1}
Let $\chi$ be a primitive ray class character of conductor $\frf$, 
and $W(\chi)$ the constant in the functional equation (\ref{eq:FE}). 
\begin{enumerate}
\item If $\bigl(\chi(\frC_1),\chi(\frC_2)\bigr)=(+1,-1)$, then 
\[L(1,\chi)=-\frac{\pi W(\chi)}{\sqrt{DN(\frf)}}
\sum_{\frC\in\Cl_K(\frf)}\chi^{-1}(\frC)\log X_1(\frC). \]
\item If $\bigl(\chi(\frC_1),\chi(\frC_2)\bigr)=(-1,+1)$, then 
\[L(1,\chi)=-\frac{\pi W(\chi)}{\sqrt{DN(\frf)}}
\sum_{\frC\in\Cl_K(\frf)}\chi^{-1}(\frC)\log X_2(\frC). \]
\end{enumerate}
\end{cor}

\begin{rem}\label{rem:Contribution}\upshape
From Theorem \ref{thm:zeta(0,C)} and Corollary \ref{cor:b_chi=1}, 
we can say that, in a sense, only infinite places for which $\chi$ is positive 
contribute to the value $L(1,\chi)$. 
We expect that there is the same principle for any totally real number field. 
\end{rem}

\subsection{An example}
We conclude our discussion with an example, 
illustrating the results in this section. 

Set $K=\Q(\sqrt{5})$ and $\frf=(4-\sqrt{5})$. 
The fundamental unit and the totally positive one are 
\[\eps_0=\frac{1+\sqrt{5}}{2},\qquad \eps=\frac{3+\sqrt{5}}{2}. \]
We also have $\eps_\frf=\eps^5$, hence $r=5$. 
Moreover, for $\frC=[O_K]$, we can take 
\begin{gather*}
\fra:=O_K,\qquad
\frb:=\base{1,\eps}=\biggl(\frac{4+\sqrt{5}}{11}\biggr)\fra^{-1}\frf,\\
\omega_0:=\eps=\mcf{3},\qquad 
z:=\frac{4+\sqrt{5}}{11}=\frac{2}{11}\omega_0+\frac{1}{11}. 
\end{gather*}
Then $\omega_k=\eps$ and $b_k=3$ for all $k\in\Z$, while 
\[(x_k,y_k)=\biggl(\frac{2}{11},\frac{1}{11}\biggr),\ 
\biggl(\frac{7}{11},\frac{9}{11}\biggr),\ 
\biggl(\frac{8}{11},\frac{4}{11}\biggr),\ 
\biggl(\frac{6}{11},\frac{3}{11}\biggr),\ 
\biggl(\frac{10}{11},\frac{5}{11}\biggr) \]
for $k\equiv 0,\,1,\,2,\,3,\,4\pmod{11}$, respectively. 
Hence $X_1(\frC)$ is defined by 
\begin{align*}
&X_1(\frC)\\
&=\Sine\biggl(\eps,\frac{7\eps+9}{11}\biggr)
\Sine\biggl(\eps,\frac{8\eps+4}{11}\biggr)
\Sine\biggl(\eps,\frac{6\eps+3}{11}\biggr)
\Sine\biggl(\eps,\frac{10\eps+5}{11}\biggr)
\Sine\biggl(\eps,\frac{2\eps+1}{11}\biggr), 
\end{align*}
and $X_2(\frC)$ is obtained by replacing $\eps$ by $\eps'$. 

Now let us exploit Theorem \ref{thm:X_i(C^*)}. 
In the present case, $\frC_1$ is the identity, since the unit $-\eps_0^5$ 
satisfies the condition of $\mu_1$: 
\[-\eps_0^5\equiv 1\mod \frf,\qquad 
-\eps_0^5<0,\qquad \bigl(-\eps_0^5\bigr)'>0. \]
Hence we have $X_2(\frC)=X_2(\frC\frC_1)^{-1}=X_2(\frC)^{-1}$, i.e.\ 
$X_2(\frC)=1$. Note that this is not trivial from its form of a product 
of the double sines. 

\begin{rem}\upshape
Shintani \cite[\S3.2]{Shintani77} proved that 
\[X(\frC)=X_1(\frC)X_2(\frC)=
\Biggl(\frac{3+\sqrt{5}}{2}-\sqrt{\frac{3\sqrt{5}-1}{2}}\Biggr)\biggm/2. \]
The fact that $X_2(\frC)=1$ seems to be new. 
\end{rem}



\end{document}